\documentstyle{amsppt}
\def\a{\alpha} \def\b{\beta} \def\d{\delta}
\def\g{\gamma}  
 \def\s{\sigma}
 \def\sm{\setminus} \def\o{\omega}
\def\to{\rightarrow}  \def\imp{\Rightarrow}
\def\r{\upharpoonright}  \def\cl{\overline}
  \def\0{\emptyset}
\def\sat{\vDash} \def\forces{\Vdash}
\def\<{\langle} \def\>{\rangle}
\topmatter

\title  Countable Toronto spaces
\endtitle

\author Gary Gruenhage and J. Tatch Moore\endauthor 

\address Department of
Mathematics, Auburn University, AL 36849\endaddress

\email garyg\@mail.auburn.edu \endemail

\address Department of Mathematics, University of Toronto, 
Toronto, Ontario M5S 1A1, CANADA \endaddress

\email justin\@math.toronto.edu \endemail


\abstract  A space $X$ is called an {\it $\alpha$-Toronto space} 
if $X$ is scattered of Cantor-Bendixson rank $\a$ and  
is homeomorphic to each of its subspaces of 
same rank. We answer a question of Steprans by constructing 
a countable $\a$-Toronto space for each $\a\leq \o$.  We also construct 
consistent examples of countable $\a$-Toronto spaces for each
$\a<\o_1$.
\endabstract
\thanks Research of the first author 
partially supported by NSF DMS-9704849 \endthanks
\endtopmatter

\document
\head 1. Introduction \endhead

For a space $X$, let $X^0=X$, let $X^{\a+1}$ be the non-isolated points
of $X^\a$, and for $\a$ a limit, let $X^\a=\bigcap_{\b<\a}X^\b$.  $X$ is 
{\it scattered} if $X^\a=\0$ for some $\a$. 
In this case, we call the least such $\a$ the {\it Cantor-Bendixson rank}, 
$r(X)$, of $X$.  The set $X^\a\sm X^{\a+1}$ (i.e., the isolated points
of $X^\a$) is called the {\it $\a^{th}$ level} of $X$.  

  The so-called Toronto problem, posed by J. Steprans \cite{S}, asks if
there is an uncountable non-discrete Hausdorff space $X$ which is
homeomorphic to each of its uncountable subspaces.  Such a space, if it exists,
is called a {\it Toronto space}. 
According to the
folklore, a Toronto space must be scattered of Cantor-Bendixson rank
$\omega_1$, and hereditarily separable (in particular, each level must
be countable). The Toronto problem is still unsettled, though it is known that
there are no Toronto spaces under $CH$, or under $PFA$ if regularity is
assumed.  

  Taking a cue from the structure of a Toronto space, Steprans calls a 
space $X$  an {\it $\a$-Toronto space} if $X$ 
is scattered of rank $\a$, and $X$ is homeomorphic to each
of its subspaces of the same rank.  For example, a convergent sequence is 
2-Toronto.  Steprans asks if 
there is an $\o$-Toronto space, and mentions that it is unknown if there
is an $\a$-Toronto space for any $\a\geq 3$.  The main result of this 
paper is 
that there is a countable $\a$-Toronto space for any
$\a\leq \o$, and there are consistent examples of countable $\a$-Toronto spaces
for each $\o<\a<\o_1$. 

Given a filter $\Cal F$ on $\o$, let $\Cal F^+$ be all $X\subset \o$
such that $X\cap F\neq \0$ for every $F\in \Cal F$.  
$\Cal F$ is said to be {\it homogeneous} if 
for any $X\in \Cal F^+$, 
the restriction of the filter $\Cal F$ to $X$ is
isomorphic to $\Cal F$.  
Let us also denote by $\Cal F\times \o$ the filter on $\o\times \o$
generated by sets of the form $F\times \o$, $F\in \Cal F$.  
We will show in Section 2 that if there is a homogeneous filter 
$\Cal F$ on $\o$ which is 
isomorphic to the filter 
$\Cal F\times \o$ on $\o^2$, then there is an $\a$-Toronto 
space for every $\a\leq \o$. If $\Cal F$ has a certain additional property 
(see Theorem 2.11), then there is an $\a$-Toronto space for every 
$\a<\o_1$.

In Section 3 we describe, in $ZFC$,  a filter $\Cal F$ which is homogeneous and 
isomorphic to $\Cal F \times \o$.  Thus there are countable  
$\alpha$-Toronto spaces in $ZFC$ for every $\alpha\leq \omega$.  This 
filter does not have the additional property required
to build countable $\alpha$-Toronto spaces for $\omega<\alpha<\omega_1$.  In
Section 4, we show that it is consistent for there to be a filter having also
the additional property, 
and hence consistent for there to be countable
$\alpha$-Toronto spaces for every $\alpha<\omega_1$.  We don't know if
such filters (i.e., having the additional property) exist in $ZFC$, or 
if there can be 
countable $\alpha$-Toronto spaces in $ZFC$ for $\omega<\alpha<\omega_1$.

Our construction also produces (consistent) 
$\a$-Toronto spaces for $\a=\o_1$ and 
$\a= \o_1+1$, albeit with uncountable levels.  We don't know if there
can be an $\o_1$-Toronto space with countable levels.

In Section 5, we show that a very natural construction of spaces 
from a filter which is somewhat different from our construction in Section 2 
cannot produce $\alpha$-Toronto spaces for $\alpha\geq 4$, though it 
does give other 3-Toronto spaces in $ZFC$.

 \head 2. Toronto spaces from the filter \endhead

Given a filter $\Cal F$ on $\o$, we will define corresponding scattered 
spaces $T_\a$, $\a$ an ordinal, of rank $\a$.  $T_\a$ will be countable 
iff $\a<\o_1$.   Most other properties of the 
spaces $T_\a$ will depend on properties of the filter $\Cal
F$.  We will show (Theorem 2.10) that if $\Cal F$ is homogeneous and 
isomorphic to the filter $\Cal F \times \o$, then $T_\a$ is an
$\a$-Toronto space for $\a\leq \o$.  If $\Cal F$ has an additional
property (see Theorem 2.11), then $T_\a$ is $\a$-Toronto for all
$\a<\o_1$.

If $(X_n)_{n\in\o}$ is a sequence of (disjoint) spaces, and 
$\Cal F$ is a filter on $\o$, let $(X_n)_{n\in\o}\cup_{\Cal F}\{\infty\}$ 
denote the space whose set is $\{\infty\}\cup\bigcup_{n\in\o}X_n$, such 
that each $X_n$ is a clopen subspace, and a neighborhood of $\infty$ has 
the form $\{\infty\}\cup \bigcup_{n\in F}X_n$, where $F\in \Cal F$.
If $\Cal F$ is a filter on a set $A$, then
$(X_a)_{a\in A}\cup_{\Cal F}\{\infty\}$ is defined similarly.

Let $T_0$ be the empty space and $T_1$ a single point space.  Let 
$T_{\a+1}=(\{n\}\times T_\a)_{n\in\o}\cup_{\Cal F}\{\infty\}$.  If $\a$ is
a limit ordinal and $T_\b$ has been defined for all $\b<\a$, let
$T_\a=\oplus_{\b<\a} T_\b$. Clearly $T_\a$ is scattered of rank $\a$,
and is countable iff $\a<\o_1$.  

To verify other properties of the spaces $T_\a$, given properties of the
filter $\Cal F$, it will be helpful to establish some general facts
about the ``$\cup_{\Cal F}$" construction.  

\proclaim{Lemma 2.1} If $W\subset (X_n)_{n\in\o}\cup_{\Cal
F}\{\infty\}$, then $$\infty\in \cl{W\sm\{\infty\}} \iff \{n:W\cap
X_n\neq\0\}\in \Cal F^+. $$  
\endproclaim
\demo{Proof} Clear. \qed\enddemo

\proclaim{Lemma 2.2} If $h_n:X_n\to Y_n$ is a homeomorphism for each 
$n\in \o$, then the map $H:(X_n)_{n\in\o}\cup_{\Cal F}\{\infty\}\to
(Y_n)_{n\in\o}\cup_{\Cal F}\{\infty\}$ defined by $H(\infty)=\infty$ and 
$H(x)=h_n(x)$ for $x\in X_n$ is a homeomorphism.
\endproclaim
\demo{Proof}.  That $H$ is a bijection and continous at every point except
(possibly) $\infty$ is clear.  It follows from Lemma 2.1 that
 $\infty\in \cl{W}\iff \infty\in \cl{H(W)}$.  Thus $H$ is also continuous 
at $\infty$.  The proof that $H^{-1}$ is continuous is similar.  Thus $H$ 
is a homeomorphism.  \qed \enddemo

\proclaim{Lemma 2.3}  For each $A\subset \o$, 
$(X_n)_{n\in\o}\cup_{\Cal F}\{\infty\}\cong 
[((X_n)_{n\in A}\cup_{\Cal F\r A}\{\infty\})\oplus((X_n)_{n\in \o\sm A}
\cup_{\Cal F\r(\o\sm A)}\{\infty\})]/\{\infty\}$.
\endproclaim    
\demo{Proof}  That the obvious mapping is a homeomorphism follows easily 
from Lemma 2.1, noting that  $B\in \Cal F^+$ iff either $B\cap A \in 
\Cal F^+ $ or $B\cap (\o\sm A)\in \Cal F^+$.  \qed 
\enddemo  

Another straightforward application of Lemma 2.1 shows:

\proclaim{Lemma 2.4}  For spaces $X_n$ and $Y_n$, $n\in\o$, 
we have \newline $(X_n\oplus Y_n)_{n\in\o}\cup_{\Cal F}\{\infty\}\cong 
[((X_n)_{n\in \o}\cup_{\Cal F}\{\infty\})\oplus((Y_n)_{n\in \o}
\cup_{\Cal F }\{\infty\})]/\{\infty\}$. 
\endproclaim

\proclaim{Lemma 2.5} If $\Cal F$ is homogeneous, then for any $A\in \Cal
F^+$, we have $$(X_n)_{n\in\o}\cup_{\Cal F}\{\infty\}\cong 
(Y_n\oplus Z_n)_{n\in \o}\cup_{\Cal F}\{\infty\}$$
where $(Y_n)_{n\in\o}$ is a re-indexing of $(X_n)_{n\in A}$, and
$Z_n=X_n$ for $n\in\o\sm A$ while $Z_n=\0$ for $n\in A$. 
\endproclaim
\demo{Proof} By homogeneity, $(X_n)_{n\in A}\cup_{\Cal F\r A}\{\infty\}$
is homeomorphic to $(Y_n)_{n\in\o}\cup_{\Cal F}\{\infty\}$
 for some re-indexing $(Y_n)_{n\in \o}$
of $(X_n)_{n\in A}$.  Also, it is clear that $(X_n)_{n\in \o\sm A}
\cup_{\Cal F\r(\o\sm A)}\{\infty\}$ is homeomorphic to $(Z_n)_{n\in \o}
\cup_{\Cal F}\{\infty\}$ where the $Z_n$'s are as given in the statement
of the lemma.  The lemma now follows by applying Lemma 2.3, the above remarks,
 and Lemma 2.4.  \qed
\enddemo

\proclaim{Lemma 2.6} If $\Cal F$ is homogenous and not the co-finite
filter, then $T_\a$ is homeomorphic to every topological sum of $T_\a$
and countably many $T_\b$ for $\b<\a$.  
\endproclaim
\demo{Proof} The result is obvious for $T_0$ and $T_1$.  Suppose it holds 
for all $\b<\a$.  If $\a$ is a limit ordinal, then $T_\a$ is by definition
the topological sum of the $T_\b$'s for $\b<\a$, and the result follows
easily from the induction hypothesis.

Suppose $\a$ is a successor, say $\a=\g+1$, and consider $T_\a\oplus 
(\oplus_{n\in\o} X_n)$, where each $X_n$ is homeomorphic to some $T_{\b_n}$, 
$\b_n<\a$.    
Note that $\Cal F$  homogenous and not the co-finite
filter implies that  $T_\a$ is homeomorphic the topological sum of itself 
and countably many  copies $Y_n$, $n\in\o$, of $T_\g$.  Now 
$T_\a\oplus (\oplus_{n\in\o} X_n)\cong [T_\a\oplus 
(\oplus_{n\in\o} Y_n)]\oplus ( \oplus_{n\in\o} X_n)\cong   T_\a\oplus 
(\oplus_{n\in\o}(X_n\oplus Y_n)$.  By the induction hypothesis, each 
$X_n\oplus Y_n$ is homeomorphic to one or two copies of $T_\g$, and the 
result follows. \qed
\enddemo

\proclaim{Lemma 2.7}  For any filter $\Cal F$ on $\o$, and disjoint spaces
$\{X_{n,m}:n,m\in\o\}$, $$(X_{n,m})_{(n,m)\in\o^2}\cup_{\Cal F\times
\o}\{\infty\}\cong (\oplus_{m\in\o}X_{n,m})_{n\in\o}\cup_{\Cal
F}\{\infty\}.$$
\endproclaim
\demo{Proof} Note that for $A\subset \o^2$, $A\in (\Cal F
\times \o)^+\iff \pi_1(A)\in \Cal F^+$.  Now it is easy 
to use Lemma 2.1 to verify that the obvious mapping is a homeomorphism. \qed
\enddemo

The next lemma shows that if in the definition of $T_{\a+1}$, each copy 
$\{n\}\times T_\a$ of $T_\a$ is replaced by the sum of finitely many or countably 
infinitely many copies of $T_\a$, the result is homeomorphic to $T_{\a+1}$,
provided $\Cal F$ has the stated properties.   
  \proclaim{Lemma 2.8} If $\Cal F$ is homogeneous and isomorphic to $\Cal
F\times\o$, then $$(\{n\}\times T_\a)_{n\in\o}\cup_{\Cal
F}\{\infty\}\cong (\oplus_{m<k_n}\{n,m\}\times T_\a)_{n\in\o}\cup_{\Cal
F}\{\infty\},$$ where $0<k_n\leq\o$.  
\endproclaim
\demo{Proof} Let $A=\bigcup_n \{n\}\times k_n$.  Then $A\in (\Cal
F\times\o)^+$.  By the assumptions on $\Cal F$,  $(\Cal
F\times\o)\r A$ is isomorphic to $\Cal F\times \o$ and $\Cal F$. 
Use an isomorphism between   $(\Cal
F\times\o)\r A$ and $\Cal F$ to constuct the natural bijection between the 
spaces.  Then use $B\in ((\Cal F\times \o)\r A)^+ \iff \pi_1(B)\in \Cal F^+$
to show via Lemma 2.1 that this bijection is a homeomorphism.
\qed
\enddemo

\proclaim{Lemma 2.9}  Suppose $\Cal F$ is homogeneous and isomorphic to
$\Cal F\times\o$.  If $A\in \Cal F^+$, then $$(X_n)_{n\in\o}\cup_{\Cal F}
\{\infty\}\cong  (\oplus_{m\in\o} Y^m_n)_{n\in \o}\cup_{\Cal F}\{\infty\},$$
where $(Y_n^m)_{n,m\in\o}$ is a reshuffling of the $X_i$'s, and for each 
$n$, \newline $|\{m:Y_n^m=X_i\text{ for some }i\in A\}|=\o$.  
\endproclaim
\demo{Proof}  We have $(X_n)_{n\in\o}\cup_{\Cal F}\{\infty\}\cong 
[((X_n)_{n\in A}\cup_{\Cal F\r A}\{\infty\})\oplus((X_n)_{n\in \o\sm A}
\cup_{\Cal F\r(\o\sm A)}\{\infty\})]/\{\infty\}\cong
[((\oplus_{m\in\o}W^m_n)_{n\in \o}\cup_{\Cal F}\{\infty\})\oplus((Z_n)_{n\in \o}
\cup_{\Cal F }\{\infty\})]/\{\infty\}$, where the $W_n^m$'s are a reindexing
of the $X_n$'s, $n\in A$, and the $Z_n$ are as is Lemma 2.5.  The first
homeomorphism exists by Lemma 2.3, and the second follows from
 $\Cal F$ being homogeneous and isomorphic to $\Cal
F\times\o$, and Lemma 2.7.  Now by Lemma 2.4, the latter quotient space is 
homeomorphic to  
$(Z_n\oplus (\oplus_{m\in\o}W^m_n)_{n\in \o}\cup_{\Cal F}\{\infty\}$.
Finally, for each $n$ let $Y_n^m$, $m\in\o$, be a reindexing of 
 $\{Z_n\}\cup\{W_n^m:m\in\o\}$. \qed
\enddemo

\proclaim{Theorem 2.10} If $\Cal F$ is homogeneous and isomorphic to
$\Cal F\times\o$, then $T_n$ is a countable $n$-Toronto space
 for every $n\leq\o$, and any subspace of $T_n$ of rank $j<n$ is 
homeomorphic to a topological sum of countably many copies of $T_j$. 
\endproclaim 
\demo{Proof} Clearly the theorem holds for $n\leq 1$.  Also, if it holds
for all $n<\o$, it is easy to check that it holds for $n=\o$, since
$T_\o$ is just the topological sum of the $T_n$'s, $n<\o$.  

We will complete the proof by showing that the theorem holds for 
$n=k+1$, given it holds for $n\leq k$, where $k<\o$. To this end,
let $X$ be a subspace of $T_{k+1}$.  If $\infty$ is either 
not in $X$ or is an isolated point of $X$, it is easy to use the
induction hypothesis to verify the conclusion of the theorem.  So,
suppose $\infty$ is a limit point of $X$.   
Let $X_n=X\cap(\{n\}\times T_k)$.  Note that $X\cong
(X_n)_{n\in\o}\cup_{\Cal F}\{\infty\}$.  Let $j_0\leq k$ be maximal such
that $\{n:r(X_n)=j_0\}\in \Cal F^+$.  Then there is $F_0\in \Cal F$ such
that $j_0=max\{r(X_n):n\in F_0\}$, and $A_0=\{n\in F_0:r(X_n)=j_0\}\in 
\Cal F^+$.  Let $N_0$ be the clopen neighborhood $\{\infty\}\cup(\bigcup_
{n\in F_0}X_n)$ of $\infty$ in $X$; note $N_0\cong 
(X_n)_{n\in F_0}\cup_{\Cal F}\{\infty\}$.    By homogeneity and Lemma 2.5,
$N_0$ is homeomorphic to 
$(Y_n\oplus Z_n)_{n\in \o}\cup_{\Cal F}\{\infty\}$
where $(Y_n)_{n\in\o}$ is a re-indexing of $(X_n)_{n\in A_0}$, and
each $Z_n$ is either $\0$ or $X_m$ for some $m\in F_0\sm A$.
By the induction hypothesis, since $r(Z_n)\leq r(Y_n) = j_0$, 
each $Y_n\oplus Z_n$ is a topological sum
of countably many copies of $T_{j_0}$.  So by Lemma 2.8, $N_0$ is
homeomorphic to $T_{j_0+1}$.  Let $N_1=X\sm N_0$.  By the induction
hypothesis, $N_1$ is homeomorphic to a topological sum of countably many
copies of $T_{j_1}$ for some $j_1\leq k$.  Of course $X\cong N_0\oplus N_1$. 
If $j_0=k$, which happens iff
$r(X)=k+1$, then $N_0$ and (by Lemma 2.6) $X$ are homeomorphic to
$T_{k+1}$.  If $j_0<k$, then by Lemma 2.6 $X$ is homeomorphic to
a topological sum of countably many copies of $T_{max\{j_0+1,j_1\}}$.
That completes the proof.\qed
\enddemo

\proclaim{Theorem 2.11} Suppose $\Cal F$ is homogeneous and isomorphic to
$\Cal F\times\o$, and satisfies 
\roster
\item"($*$)" Whenever $f:\o\to\o$ is unbounded on every $F\in \Cal F$,
there is some $A\in \Cal F^+$ such that $f\r A$ is finite-to-one.
\endroster
Then $T_\a$ is a countable $\a$-Toronto space for every $\a<\o_1$, and
any subspace of $T_\a$ of rank $\b<\a$ is homeomorphic to a 
topological sum of countably many copies of $T_\b$.  
\endproclaim
\demo{Proof} Theorem 2.10 shows this theorem holds for $\a\leq\o$.
Suppose it holds for all $\b<\a$, where $\a<\o_1$, 
and consider $X\subset T_\a$.

If $\a$ is a limit, then $T_\a=\oplus_{\b<\a} T_\b$. 
 By the induction
hypothesis, each subspace $X\cap T_\b$ of $X$ is homeomorphic to the
topological sum of copies of $T_{\b'}$ for some $\b'<\b$.  Now one can
use Lemma 2.6 to split up or group these $T_{\b'}$'s in the appropriate
way to show that $X\cong T_{r(X)}$ if $r(X)$ is a limit ordinal, and is 
homeomorphic to the sum of countably many copies of $T_{r(X)}$ otherwise.

It remains to verify the theorem in case $\a$ is a successor, say
$\a=\g+1$.  As in the proof of the Theorem 2.10, we may suppose that
$\infty$ is a limit point of $X$, in which case $X\cong
(X_n)_{n\in\o}\cup_{\Cal F}\{\infty\}$, where 
$X_n=X\cap(\{n\}\times T_\g)$. 

Let $\d_n=r(X_n)$, and let $\d$ be
minimal such that, for some $F_0\in \Cal F$, $$sup\{\d_n:n\in F_0\}=
\d.$$  Let $A=\{n\in F_0:\d_n=\d\}$.  If $A\in \Cal F^+$, let $A_0=A$.
If $A\not\in\Cal F^+$, then we may assume $A=\0$, i.e., $\d_n<\d$ for 
each $n\in F_0$.  Note that therefore $\d$ is a limit ordinal.  Let $\b_n$,   
$n\in \o$, be increasing with supremum $\d$.  Define $f:\o\to\o$ such
that $f(n)=m$ iff $\b_m\leq\d_n<\b_{m+1}$.  By minimality of $\d$, $f$ is
unbounded on every $F\in \Cal F$. By ($*$), there is $A_0\in\Cal F^+$, 
$A_0\subset F_0$, such that $f\r A_0$ is finite-to-one.  

Thus, whether $A\in \Cal F^+$ or not, we have an $A_0\in\Cal F^+$, 
$A_0\subset F_0$, such
that $sup\{\d_n:n\in B\}=\d$ for every infinite $B\subset A_0$.  
 Let $N_0$ be the clopen neighborhood 
$\{\infty\}\cup(\bigcup_
{n\in F_0}X_n)$ of $\infty$ in $X$.
By homogeneity and Lemma 2.9, $N_0\cong (Z_n)_{n\in\o}\cup_{\Cal F}
\{\infty\}$, where each $Z_n$ is a topological sum of $X_m$'s, $m\in F_0$,
with $m\in A_0$ infinitely often.
 Thus $r(Z_n)=\d$ for all $n$.  By the induction hypothesis, $Z_n$ is 
homeomorphic to a topological sum of countably many copies of $T_\d$, 
hence by Lemma 2.8, $N_0\cong T_{\d+1}$.  The proof is now completed as 
in Theorem 2.10.\qed
\enddemo
  
We now show that the condition ($*$) is necessary for $T_{\o+1}$ to
be $(\o+1)$-Toronto.  

\proclaim{Theorem 2.12}  Suppose $\Cal F$ is a filter on $\o$ which
fails to satisfy condition ($*$) of Thereom 2.11.  Then $T_{\o+1}$ is
not $(\o+1)$-Toronto.  
\endproclaim
\demo{Proof}   Recall $T_{\o+1}=(\{n\}\times
T_\o)_{n\in\o}\bigcup_{\Cal F}\{\infty\}$.  Suppose $T_{\o+1}$ is
$(\o+1)$-Toronto, and let $f:\o\to\o$ be such that $|f(F)|=\o$ for
every $F\in \Cal F$.  We will show that $f\r A$ is finite-to-one
for some $A\in \Cal F^+$. 
To this end, let  $X=(\{n\}\times T_{f(n)})_{n\in\o}\bigcup_{\Cal
F}\{\infty\}$.  Then $X$ is homeomorphic to a subspace of
$T_{\o+1}$, and it follows from the fact that  $|f(F)|=\o$ for
every $F\in \Cal F$ that $r(X)=\o+1$.  

 So there must exist a homeomorphism $h:T_{\o+1}\to X$.  For $k\in
\o$ let $$B_k=\{n\in\o: h(\{k\}\times T_\o)\cap (\{n\}\times
T_{f(n)})\neq \0.$$  Since $r(T_\o)=\o$, we must have
$|f(B_k)|=\o$.  Thus we can choose $n_k\in B_k$ such that
$f(n_0)<f(n_1)<...$.  Now pick $x_k\in      h(\{k\}\times T_\o)\cap
(\{n_k\}\times T_{f(n_k)})$.  Note that $\infty \in cl(h^{-
1}(\{x_k:k\in\o\}))$.  It follows that $A=\{n_k:k\in \o\}$ is in
$\Cal F^+$.  Since $f$ is one-to-one on $A$, that completes the
proof. \qed

\enddemo

\proclaim{Theorem 2.13} If $\Cal F$ satisfies the conditions of 
Theorem 2.11, then $T_{\o_1}$ and  $T_{\o_1+1}$ are $\o_1$-Toronto and
$(\o_1+1)$-Toronto, respectively (albeit with uncountable levels).
However, $T_{\o_1+2}$ is not $(\o_1+2)$-Toronto.  
\endproclaim
\demo{Proof} That $T_{\o_1}$  is $\o_1$-Toronto follows just like the
limit case of the proof of the preceding theorem.  Now consider a
subspace $X$ of $T_{\o_1+1}$ of rank $\o_1+1$, and let $X_n=X\cap 
(\{n\}\times T_{\o_1})$.  Then for $\Cal F^+$-many $n$'s, $r(X_n)=\o_1$
and for such $n$, $X_n\cong T_{\o_1}$.  Now use Lemmas 2.5 and 2.8 
 as in Theorem 2.10 to complete the proof that $X\cong T_{\o_1+1}$.  
Thus $T_{\o_1+1}$ is $(\o_1+1)$-Toronto.  

Let $Y$ be the subspace of $T_{\o_1+1}$ consisting of only its 
isolated points and the point $\infty$.  Note that every neighborhood of 
$\infty$ in $Y$ is uncountable.  It follows that  $T_{\o_1+2}$ has a 
subspace $Z$ of rank $\o_1+2$ with a point at level 1 every neighborhood of 
which is uncountable.  But every level 1 point of $T_{\o_1+2}$ has a   
countable neighborhood.  Thus $Z\not\cong T_{\o_1+2}$ and so $T_{\o_1+2}$ 
is not $(\o_1+2)$-Toronto. \qed
\enddemo

\head 3. The ZFC filter  \endhead

In this section we will define a homogeneous filter $\Cal F$ on
$\omega$ such that $\Cal F$ is isomorphic to $\Cal F \times \omega$.
It will be more convenient for our purposes to actually define
the filter on the countable ordinal $\omega^\omega = \sum_n \omega^n$.
The filter $\Cal F$ is the collection of all $A \subseteq \omega^\omega$
such that the order type of $\omega^\omega \setminus A$ is less
than $\omega^\omega$.
Note that the $\Cal F$-positive subsets of $\omega^\omega$ 
are simply those
which have order type $\omega^\omega$.

\proclaim{Theorem 3.1}
$\Cal F$ is homogeneous and isomorphic to $\Cal F\times \omega$.
\endproclaim

\demo{Proof}
Notice that
if $X \subseteq \omega^\omega$ is in  $\Cal F^+$ and
$\Phi:X \to \omega^\omega$ is an order isomorphism then
$\Phi$ is also an isomorphism between $\Cal F \restriction X$
and $\Cal F$.  Thus $\Cal F$ is homogenous.

To see that $\Cal F$ is isomorphic to $\Cal F\times \omega$, recall 
that the ordinal $\omega \cdot \omega^\omega$ is
the order type of the set $\omega \times \omega^\omega$ equiped with the
lexicographical order.
It is easy to see that
$$\omega \cdot \omega^\omega = \sum_{n=1}^\infty
\omega \cdot \omega^n = \omega^\omega.$$
Let $\Phi$ be an order isomorphism between
$\omega^\omega \times \omega$ with the reverse lexicographical order
and $\omega^\omega$.
Define $\Cal F^*$ to be the preimage of $\Cal F$ under $\Phi$.
It now suffices to show that $\Cal F^*$ is equal to $\Cal F \times \omega$.
It should be clear that $\Cal F^*$ contains $\Cal F \times \omega$.
If $A$ is in $\Cal F^*$ then define $A^0$ to be the union of all 
$E \subseteq A$ such
that $E = \pi^{-1}(\pi(E))$
where $\pi:\omega^\omega \times \omega \to \omega^\omega$ is the
projection onto the first coordinate.
Notice that $A^0$ is also in $\Cal F$ since if the complement of $A$ has
order type $\alpha < \omega^\omega$ then the complement of $A^0$
has order type at most $\omega \cdot \alpha < \omega^\omega$.
Since $\pi(A^0)$ must be in $\Cal F$, and since $\pi^{-1}(\pi(A^0)) =  A^0$,
$A \supseteq A^0$ is in $\Cal F \times \omega$ and we are finished. \qed
\enddemo

\proclaim{Corollary 3.2}  There are, in $ZFC$, $n$-Toronto spaces for
every $n\leq \o$. 
\endproclaim
\demo{Proof}  Immediate from Theorems 3.1 and 2.10.
\enddemo

Unfortunately, this filter does not produce $\a$-Toronto 
spaces for $\a>\o$ by Theorem 2.12 and:

\proclaim{ Proposition 3.3} Let $\Cal F$ be the filter of Theorem 3.1.  
Then $\Cal F$ does not satisfy condition ($*$) of Theorem 2.11. 
\endproclaim
\demo{Proof}  Define $f:\o^\o\to \o$ by $f(\a)=k$ iff $\a\in
[\o^k,\o^{k+1})$.  If the restriction of $f$ to a set $A\subset \o^\o$ is
finite-to-one, clearly $A$ has order type $\o$ and hence is not in $\Cal
F^+$.  \qed
\enddemo

In the next section, we will show that at least there are consistent
examples of filters satisfying all the conditions of Theorem 2.11. 

\head 4. The consistent filter \endhead 

  The purpose of this section is to prove the following.

\proclaim{Theorem 4.1}  If ZFC is consistent, then it is consistent with 
ZFC that there is a filter $\Cal F$ on $\o$ satisfying:
\roster
\item"(i)" $\Cal F$ is homogeneous;
\item"(ii)" $\Cal F$ is isomorphic to the filter $\Cal F\times\o$ on $\o^2$;
\item"(iii)" Whenever $f:\o\to \o$ is unbounded on every $F\in \Cal F$, 
then there is some $A\in \Cal F^+$ such that $f\r A$ is finite-to-one.
\endroster
\endproclaim

  The  starting point for our construction is the following observation.  
Suppose $\Cal F_e=\{F_{e\a}:\a<\o_1\}$, 
$e=0,1$, are subbases for filters on $\o$ such that, 
whenever $H$ and $K$ are 
disjoint finite subsets of $\o_1$ and $ \cap_{\a\in H}F_{e\a}\sm 
\cup_{\b\in K} F_{e\b}$ is non-empty, then it is infinite, and so is the 
corresponding set using $\Cal F_{1-e}$.  Then there  is a natural 
$\sigma$-centered poset $P$ forcing
$\Cal F_0$ and $\Cal F_1$ to be isomorphic: namely, $P$ consists of all 
pairs $p=(\tau^p,H^p)$, where $\tau^p$ is a finite one-to-one function 
from $\o$ to $\o$, and $H^p\in [\o_1]^{<\o}$.  Declare $q\leq p$ if
$\tau^q
\supset \tau^p$, $H^q\supset H^p$, and for each $\b\in H^p$, we have 
$n\in F_{0 \b} \iff\tau^q(n)\in F_{1 \b}$ 
whenever $n\in dom(\tau^q\sm\tau^p)$.  This forcing adds a function
$t:\o\to\o$ such that $t(F_{0 \b})=^* F_{1 \b}$ for every $\b<\o_1$,
and it is easy to see that such a function $t$ is an isomorphism.

  Call a pair of filters having subbases as above a {\it good} pair.  
Any pair of filters generated by $\o_1$-sized independent families is 
a good pair.   Our naive idea to start with a filter $\Cal F$ 
on $\o$ generated 
by an independent family $\Cal A=\{A_\a:\a<\o_1\}$,
 consisting of $\o_1$-many Cohen reals, 
then force it to be 
isomorphic to the filter $\Cal F\times\o$  on 
$\o^2$, which is generated by the independent family
$\{A\times \o: A\in \Cal A\}$, and finally iterate the 
type of poset described in the previous paragraph $\o_1$ times (we start 
with a model of CH) to force $\Cal F$ to be homogeneous.   

  Note that for any infinite subset $X$ of $\o$ in the ground model,
the restriction of $\Cal A$ to $X$ is still an independent family.  
The problem one runs into, however,  is that this is  not true for
many subsets $X$ added by the forcings.  For example, if some infinite set
$X$ is added which is almost contained in every member of $\Cal A$, then 
the restriction of $\Cal F$ to $X$ is the cofinite filter,
and then there is no hope of making $\Cal F$ homogeneous.  So we must 
in particular show 
sets like this are not added.   In fact, we show that for any subset $X$ 
of $\o$ added 
at some stage of the iteration,  if $X\in \Cal F^+$, then there is some 
$\delta<\o_1$ such that the restriction of $\{A_\a:\a\geq\delta\}$ to 
$X$ is
an independent family.  By Lemma 4.2 below, this 
turns out to be enough for 
$\Cal F$ and its restriction to $X$ to be a good pair of filters
(witnessed by the subbase $\{\bigcap_{\phi(\a)<n}A_\a\}_{n<\o}\cup
\{A_\a:\a\geq\d\}$ for $\Cal F$ and its restriction to $X$ for 
$\Cal F\r X$),
and so they can be forced to be isomorphic.  

  Forcing notation follows Kunen \cite{Ku}; in particular, $Fn(X,Y)$ denotes
the set of all functions from a finite subset of $X$ into $Y$.  
For $A\subset \o$, 
we let $A^1=A$ and $A^0=\o\sm A$.

\proclaim{Lemma 4.2}  Let $\Cal A=\{A_\a:\a<\o_1\}$ be an independent 
family of subsets of $\o$, and let $\Cal F$ be the filter generated by 
$\Cal A$.  Given $\rho\in Fn(\o_1,2)$, let $L_\rho=\bigcap_{\a\in dom(\rho)}
A_\a^{\rho(\a)}$.  Suppose the following holds: 
$$\forall X\subset \o[X\in \Cal F^+\imp \exists \gamma<\o_1(X\cap L_\rho
\neq\0 \text{ for all } \rho\in Fn(\o_1\sm\gamma,2))].$$ 
(To express the property in words, one might say it means that 
the restriction of 
$\Cal A$ to any member of $\Cal F^+$ is ``eventually independent".)

Then for every $X\in \Cal F^+$, there exists $\delta<\o_1$ and a finite-to-one
$\phi:\delta\to \o$ such that: $$[\bigcap_{\phi(\a)<n}A_\a\sm  
\bigcap_{\phi(\a)\leq n}A_\a]\cap L_\rho\cap X\neq \0$$ for all $n<\o$ and 
$\rho\in Fn(\o_1\sm \delta,2)$.  
\endproclaim

\demo{Proof}  Let $X\in \Cal F^+$.
 Let $M$ be a countable elementary submodel containing $\Cal A, X$, and a
function $Y\mapsto \gamma(Y)\in \o_1$ witnessing the hypothesized
property.  Let $\delta=M\cap \o_1$.  Construct a finite-to-one fuction
$\phi:\delta\to \o$ such that, for each $n$, $\phi^{-1}(n)\not\subset 
\gamma(X\cap \bigcap_{\phi(\a)<n}A_\a)$ as follows.  Let
$\d=\{\d_0,\d_1,...\}$.  Note that $A_\a\in M$ for each $\a<\d$, and
$\g(Y)<\d$ for each $Y\in M\cap \Cal F^+$.  Hence we can inductively
choose a finite subset $\phi^{-1}(n)$ of $\d$ containing:
\roster
\item"(1)" $\d_i$, where $i$ is least such that $\d_i\not\in
\bigcup_{j<n}\phi^{-1}(j)$;
\item"(2)" $\d_k>\g(X\cap\bigcap_{\phi(\a)<n} A_\a)$.
\endroster

We claim that this $\delta$ and $\phi$ have the desired properties.
To see this, fix $n\in\o$.  Let $Z_n= [\bigcap_{\phi(\a)<n}A_\a\sm  
\bigcap_{\phi(\a)\leq n}A_\a]\cap X$, and let 
$\rho\in Fn(\o_1\sm\delta,2)$.  
We need to show $Z_n\cap L_\rho\neq \0$.  
Choose $\a'\in \phi^{-1}(n)\sm \gamma(X\cap \bigcap_{\phi(\a)<n}A_\a)$. Let 
$\rho'=\rho^\frown \langle\a',0\rangle $.  Since $\rho'\in Fn(\o_1\sm 
\gamma(X\cap \bigcap_{\phi(\a)<n}A_\a))$, we have 
 $\0\neq [X\cap\bigcap_{\phi(\a)<n}
A_\a]\cap L_{\rho'}=  [(X\cap\bigcap_{\phi(\a)<n}
A_\a)\sm A_{\a'}]\cap L_{\rho}\subset[(X\cap\bigcap_{\phi(\a)<n}
A_\a)\sm\bigcap_{\phi(\a)\leq n}
A_\a]\cap L_{\rho}=Z_n\cap L_\rho$.  \qed
\enddemo

Now we describe the posets that will be used in the iteration.  
Let $\Cal F$ be the filter generated by an independent family 
$\Cal A=\{A_\a:\a<\o_1\}$,
let $X\in \Cal F^+$, let $\delta<\o_1$, and let $\phi:\delta\to \omega$
be a finite-to-one function satisfying the conclusion of Lemma 4.2.
Let $Q(\Cal A,X,\delta,\phi)$ be the poset consisting of all
$p=\langle\tau^p,F^p,n^p\rangle $ such that:
\roster
\item"(a)" $\tau^p$ is a finite one-to-one function from $\o$ to $X$;
\item"(b)" $F^p\in [\o_1\sm \delta]^{<\o}$;
\item"(c)" $n^p\in\o$.
\endroster
Define $q\leq p$ iff $\tau^q\supset\tau^p$, $F^q\supset F^p$, $n^q\geq
n^p$, and
\roster
\item"(i)"$\forall n\in dom(\tau^q\sm\tau^p)\forall \b\in F^p[n\in
A_\b\iff \tau^q(n)\in A_\b];$
\item"(ii)" $\forall n\in dom(\tau^q\sm\tau^p)
\forall k\leq n^p[n\in \bigcap_{\phi (\a)<k}A_\a
\iff \tau^q(n)\in \bigcap_{\phi (\a)<k}A_\a]$.
\endroster

\proclaim{Lemma 4.3} The poset $Q=Q(\Cal A,X,\delta,\phi)$ is $\sigma$-centered,
and if $G$ is a $Q$-generic filter, then $t=\cup\{\tau^p:p\in G\}$ is a
bijection from $\o$ to $X$ such that $t(A_\b)=^*A_\b\cap X$ for all
$\b\geq \delta$, and $t(\bigcap_{\phi (\a)<n}A_\a)=^* X\cap
\bigcap_{\phi (\a)<n}A_\a$ for all $n<\o$.  In particular, $t$
witnesses that in $V[G]$,  $\Cal F$ is isomorphic to its restriction to
$X$.  
\endproclaim
\demo{Proof}
Clearly any two conditions $p$ and $q$ for which $\tau^p=\tau^q$ are 
compatible, and so the poset is $\sigma$-centered.  

Let $G$ be a $Q$-generic filter.  First we show that, for each 
$k\in\o$, the subset of $Q$ consisting of
all $p$ with $k\in dom(\tau^p)$ is dense.
 To this end, suppose $k\not\in dom(\tau^p)$.  Let $\rho:F^p\to 2$ be
such that $k\in A_\b\iff \rho(\b)=1$.  Let $n\leq n^p$ be maximal
such that $k\in \bigcap_{\phi(\a)<n}A_\a$.  By the property of $\phi$, 
the set  $[\bigcap_{\phi(\a)<n}A_\a\sm  
\bigcap_{\phi(\a)\leq n}A_\a]\cap L_\rho\cap X$ infinite, so we can
choose a natural number $k'$ in this set which is not in $ran(\tau^p)$.  
Let $\tau^q=\tau^p\cup\{\langle k,k'\rangle \}$, $F^q=F^p$, and $n^q=n^p$.  Then 
$q\leq p$.  

  By a similar argument, for each $k\in X$,
the set of all 
$p$ with $k\in ran(\tau^p)$ is dense.  It follows that $t:\o\to X$ is a bijection. 

We now prove that $t$ is an isomorphism between $\Cal F$ and $\Cal F\r X$
by showing that 
\roster
\item"(1)" $t(A_\b)=^* X\cap A_\b$ for every $\b\geq\delta$; and
\item"(2)" $t(\bigcap_{\phi(\a)<n}A_\a)=^* X\cap\bigcap_{\phi(\a)<n} 
A_a$ for each $n<\o$.
\endroster

To see (1), fix $\b\geq \delta$.  There is $p\in G$ with $\b\in F^p$.  Let
$k\not\in dom(\tau^p)$.  There is $q\in G$ with $q\leq p$ and 
$k\in \tau^q$.  Then $k\in A_\b\iff \tau^q(k)\in A_\b\iff t(k)\in A_\b$.
Similarly, if $k\in X\sm ran(\tau^p)$, then $k\in A_\b\iff t^{-1}(k)
\in A_\b$.  It follows that
$t(A_\b\sm dom(\tau^p))=X\cap A_\b\sm ran(\tau^p)$.  
   
 Now for (2), fix $n<\o$.There is $p\in G$ with $n\leq n^p$.  Let
$k\not\in dom(\tau^p)$.  There is $q\in G$ with $q\leq p$ and 
$k\in \tau^q$.  Then $k\in \bigcap_{\phi(\a)<n}A_\a\iff 
\tau^q(k)\in \bigcap_{\phi(\a)<n}A_a\iff t(k)\in \bigcap_{\phi(\a)<n}A_a$.
The analogous statement is true for $k\in X\sm ran(\tau^p)$.  
  It follows that
$t(\bigcap_{\phi(\a)<n}A_\a\sm dom(\tau^p))=X\cap\bigcap_{\phi(\a)<n} 
A_a\sm ran(\tau^p)$.

It follows that the bijections $t$ and $t^{-1}$ map elements of 
$\Cal F$ to elements of its restriction to $X$, and vice-versa.  
So $t$ witnesses 
that these two filters are isomorphic.  \qed
\enddemo

Now we define the poset $Q_1(\Cal A)$ forcing $\Cal F$ to be isomorphic 
to the filter $\Cal F\times \o$ on $\o^2$.
Let $p=\langle\tau^p,F^p\rangle $ be in $Q_1(\Cal A)$ iff:
\roster
\item"(a)" $\tau^p$ is a finite one-to-one function from $\o$ to $\o^2$;
\item"(b)" $F^p\in [\o_1]^{<\o}$.
\endroster
Define $q\leq p$ iff $\tau^q\supset\tau^p$, $F^q\supset F^p$, 
and $\forall n\in dom(\tau^q\sm\tau^p)\forall \b\in F^p[n\in
A_\b\iff \tau^q(n)\in A_\b\times\o].$

\proclaim{Lemma 4.4} The poset $Q_1(\Cal A)$ is $\sigma$-centered,
and if $G$ is a $Q_1(\Cal A)$-generic filter, then 
$t=\cup\{\tau^p:p\in G\}$ is a
bijection from $\o$ to $\o^2$ such that $t(A_\b)=^*A_\b\times\o$ for all
$\b\in\o_1$.  In particular, $t$
witnesses that in $V[G]$,  $\Cal F$ is isomorphic to $\Cal F\times\o$.
\endproclaim
\demo{Proof}  Note that the collection $\{A_\a\times \o:\a<\o_1\}$, which 
generates the filter on $\o^2$,  is an
independent family.  Thus Lemma 4.4 follows by a proof 
similar to (and somewhat shorter than, since the complication of the
finite-to-one function is not involved here) that of Lemma 4.3.  \qed
\enddemo

Now we describe the iteration, which is a finite support iteration 
$P_{\o_1}$ of 
\newline $\langle\dot{Q_\a}\rangle _{\a<\o_1}$.  
Let the ground model $V$ satisfy $CH$.
 Let $Q_0=Fn(\o_1\times\o,2)$; i.e.,
$Q_0$ is the poset for adding $\o_1$-many Cohen reals.  If $G_0$ is
$Q_0$-generic and $\a<\o_1$, let $A_\a=\{n\in\o:\cup G_0(\a,n)=1\}$.
Then $\Cal A=\{A_\a:\a<\o_1\}$ is an independent family in $V[G_0]$.
Let $\dot{Q_1}$ be a $Q_0$-name for the forcing $Q_1(\Cal A)$ of Lemma 4.4.  
Each $\dot{Q_\a}$ for $\a>1$ will be a name for a forcing $Q(\Cal
A,X_\a,\delta_\a, \phi_\a)$ as in Lemma 4.3 forcing $\Cal F$ and its 
restriction to $X_\a$ to be isomorphic.  Since $V\sat CH$, and each
poset in the iteration is $CCC$ and has size $\o_1$, it follows that the 
final model satisfies $CH$, and we can arrange the  
$X_\a$'s to include (names for)
every $X\in \Cal F^+$ in the final model.  Thus in the end $\Cal F$ will be 
homogeneous, and, thanks to $\dot{Q_1}$, will be isomorphic to 
$\Cal F\times\o$.  It turns out nothing more need be done to obtain the 
additional property (iii) of Theorem 4.1.

However, a problem with the above simple outline is that, 
given that $X_\a\in \Cal F^+$ in $V^{P_\a}$, in order to continue the
iteration we must show that there is forced to be some $\delta_\a<\o_1$ 
and finite-to-one function $\phi_\a:\delta_\a\to \o$ such that 
$Q(\Cal A, X_\a, \delta_\a,\phi_\a)$ exists, 
i.e, the conclusion of Lemma 4.2
is satisfied.   So we must show that $P_\a$ forces the hypothesis of 
Lemma 4.2.  

To establish notation, let us describe the iteration more precisely. We 
will often abuse notation by letting sets be names for themselves, when 
this should cause no confusion, but we will also use $\dot{X}$ to denote
a name for $X$ when we want to emphasize that we are talking about
names.  

First note that, w.o.l.o.g.,  we can think of members of the iteration as 
having the form
$$p=\langle\sigma^p,\langle\tau_1^p,F_1^p\rangle ,\langle\tau_\g^p,F_\g^p,n_\g^p\rangle _{\g\in D^p}\rangle $$
where 
\roster
\item"(i)" $\s^p\in Fn(\o_1\times \o,2)$;
\item"(ii)" $\tau_1^p$ is a finite one-to-one function from $\o$ to
$\o^2$ and $F_1^p\in [\o_1]^{<\o}$;
\item"(iii)" $D^p\in [\o_1\sm 2]^{<\o}$;
\item"(iv)" For each $\g\in D^p$, $\tau_\g^p$ is a finite one-to-one function
from $\o$ to $\o$, $F_\g^p\in [\o_1]^{<\o}$, and $n_\g^p<\o$.
\endroster
 Further, there is a sequence $(\dot{X_\a},\dot{\d_\a},\dot{\phi_\a})$,
$1<\a<\o_1$, of names such that:
\roster
\item"(v)"$\forces_\a ``\dot{X_\a}\in \Cal F^+, \dot{\d_\a}<\o_1, \dot{\phi_\a}:
\dot{\d_\a}
\to\o$ is finite-to-one, and $\dot{X_\a},\dot{\d_\a},\dot{\phi_\a}$
satisfy the conclusion of Lemma 4.2",
\endroster
and for each $\g\in D^p$, 
\roster
\item"(vi)"
$p\r \g\forces ``p(\g)\in 
Q(\dot{\Cal A},\dot{X_\a},\dot{\d},\dot{\phi})".$ 
\endroster
In particular this 
implies 
\roster
\item"(vii)"$p\r\g\forces ``ran(\tau_\g^p)\subset \dot{X_\a}\text{ and } 
F_\g^p \subset \o_1\sm \dot{\d}_\a"$. 
\endroster
 As indicated above, we can also assume that
\roster
\item"(viii)"
$\forces_{\o_1}``X\in \Cal F^+\imp \exists \a<\o_1(X=\dot{X}_\a)".$    
\endroster

Let us also make the following simple but useful observation.  If
we change $p$ by extending $\s^p$ and doing nothing else, we get a 
stronger condintion; it follows that $p\forces n\in A_\b\iff \s^p\forces
n\in A_\b\iff \langle(\b,n),1\rangle \in \s^p$.   

In order to show that an iteration as described above actually exists, 
we need to show that in $V^{P_\a}$,
any $X\in \Cal F^+$ satisfies the
hypothesis of Lemma 4.2.  Then if at stage $\a$ we are given 
some $X_\a\in \Cal F^+$, we can conclude that there is a
$\d_\a<\o_1$ and finite-to-one $\phi_\a:\d_\a\to \o$ such that $X_\a,
\d_\a$, and $\phi_\a$ satisfy the conclusion of Lemma 4.2, and thus 
continue the iteration as described.

First we show, roughly speaking, that in $V^{P_\a}$, 
members of $\Cal F^+$   
are not contained in an ``orbit" of finitely many of the previously added   
isomorphisms $t_\b:\o\to X_\b$, $1<\b<\a$, and $t_1:\o\to\o^2$.  
This will be needed later to show that we are free enough to extend
conditions to force things we want to force.
 
Let $T$ be a collection of one-to-one functions 
from $\o$ to $\o$ or to $\o^2$, 
and let $k\in\o$.  We say that $O(T,k)\subset \o$ is the {\it orbit of 
$k$ under $T$} if:
\roster
\item"(i)" Whenever $t\in T$, $t:\o\to\o$, and $n\in O(T,k)$, then 
$t(n)$ and $t^{-1}(n)$ (if defined) are in $O(T,k)$;
\item"(ii)" Whenever $t\in T$, $t:\o\to\o^2$, and $n\in O(T,k)$, then 
$\pi_1(t(n))\in O(T,k)$ and $t^{-1}(n,j)\in O(T,k)$ for all $j$.  
\item"(iii)" $O(T,k)$ is the smallest set containing $k$ and satisfying
conditions (i) and (ii) above. 
\endroster
Note that $n\in O(T,k)$ iff there is a finite sequence $n_0,n_1,...,n_l$ 
of natural numbers such that $n_0=k$ and $n_l=n$, and a finite sequence 
$t_0,t_1,...,t_l$ 
of members of $T$ such that, for each $i<l$, either $n_{i+1}=t_i(n_i)$,        
 $n_{i}=t_{i+1}(n_{i+1})$, $n_{i+1}=\pi_1(t_i(n_i))$, or 
$n_{i}=\pi_1(t_{i+1}(n_{i+1})$.  

\proclaim{Lemma 4.5}  Let $t_1:\o\to\o^2$ be the bijection added by the  
first coordinate (i.e., by $\dot{Q_1}$) of the forcing $P_\a$, and 
for $1<\b$ let $t_\b:\o\to\o$ be that added by the $\b^{th}$ coordinate.       
Let $T$ be a finite subset of $\{t_\b:0<\b<\a\}$.  Then for each 
$k\in\o$, $$\forces_\a ``O(T,k)\not\in \Cal F^+".$$
\endproclaim

\demo{Proof}  Suppose not.  Then for some $p\in P_\a$, $p\forces O(T,k)\in 
\Cal F^+.$  We assume $p$ has the form described above, and we may also 
assume $p\forces T=\{t_\g:\g\in E\}$, where $E\in [\a]^{<\o}$. 
 Choose $\mu<\o_1$ 
such that no $(\mu,n)$ is in $dom(\s^p)$ and 
$p\r\g\forces \dot{\d_\g}\leq\mu$ for each $\g\in D^p$ (such $\mu$ exists by 
$CCC$).  Let $k'<\o$ be greater than $k$ and any integer mentioned in  
$\s^p$ or the $\tau_\g^p$'s, $\g\in \{1\}\cup D^p$.  Let $p'$ be
obtained from $p$ by adding $\langle(\mu,i),0\rangle $ to $\s^p$ for each $i<k'$,
and adding $\mu$ to each $F_\g^p$, $\g\in \{1\}\cup D^p$.  It is easy to
check that $p'$ is a condition stronger than $p$.  

Since $p'\forces O(T,k)\in \Cal F^+$, there are $q\leq p'$ and $m\geq k'$ 
such that $q\forces m\in O(T,k)\cap A_\mu$. By extending $q$ if 
necessary, we may assume that there is a sequence $\langle n_0,n_1,...,n_l\rangle $     
of integers and a sequence $\langle\g_0,\g_1,...\g_l\rangle \in E^{<\o}$ such that 
this sequence of integers and the sequence of members of $T$ corresponding 
to the $\g_i$'s witness that $m\in O(T,k)$ (i.e., $n_0=k$, $n_l=m$, etc...
see the discussion of $O(T,k)$   
prior to the statement of this lemma).  We may also assume that $q$,
hence $\s^q$, decides whether or not $n_i\in A_\mu$, for all $i\leq l$, 
and if, e.g., $q\forces t_{\g_i}(n_i)=n_{i+1}$, then $\tau_{\g_i}^q
(n_i)=n_{i+1}$.  
The following claim contradicts $q\forces m\in A_\mu$, proving the lemma.

{\it Claim: For each $i\leq l$, $q\forces n_i\not\in A_\mu$.}
  
{\it Proof of Claim.}  Since $\s^{p'}(\langle\mu,k\rangle )=0$, $n_0=k$, and 
$q\leq p'$, we have   $q\forces n_0\not\in A_\mu$.   
Suppose the claim is false, and let $i<l$ be least such that $q\forces 
n_{i+1}\in A_\mu$.  Suppose, e.g., that, $\g_i=1$ and  $\pi_1
(\tau_1^q(n_i))=n_{i+1}$ (all other cases, which we omit, are similar).      
Note that $q\forces 
n_{i+1}\in A_\mu$ and $q\leq p'$ implies that $n_{i+1}$ is not mentioned 
by $\tau_1^p$, and so not by $\tau_1^{p'}$ either.  Thus $n_i\not\in 
dom(\tau_1^{p'})$.  But now $q\leq p'$ 
and $\mu\in F_1^{p'}$ imply $q\r 1\forces n_i\in A_\mu\iff \tau_1^q(n_i)
\in A_\mu\times\o\iff n_{i+1}\in A_\mu$.  But by our choice of $i$,
$q\forces ``n_i\not\in A_\mu \text{ and } n_{i+1}\in A_\mu"$, contradiction.
\qed
    
\enddemo

The following lemma shows that, under certain circumstances and for certain
values of $\b$ and $m$, we are free to extend 
two compatible conditions so that $m$ is forced to be in $A_\b$, or not, 
as we wish.  

\proclaim{Lemma 4.6}  Let $\a<\omega_1$, let $q\in P_\a$, and let $M$ be 
a countable elementary submodel containing $q$ and $P_\a$.  Let $p\in P_\a$ 
be compatible with $q$ (we don't require $p\in M$ however).    
Let $l$ be greater than any integer mentioned by $p$ (i.e., greater
than any $n_\mu^p$, $\mu\in D^p$, or anything in the domain or range of 
$\s^p$ or of any $\tau_\mu^p$).  Let $T=\{t_\mu:\mu\in D^p\cup\{1\}\}$, 
where $t_\mu$ is the isomorphism added by the $\mu^{th}$ coordinate of the 
forcing. 
Suppose the following hold:
\roster
\item"(i)" $D^q\supset D^p$;
\item"(ii)" For each $\mu\in D^p\cup\{1\}$, $\tau_\mu^q\supset \tau_\mu^p$ 
and $F_\mu^q\supset F_\mu^p\cap M$;
\item"(iii)" For each $\mu\in D^p$, $n_\mu^q\geq n_\mu^p$;
\item"(iv)" For each $\mu\in D^p$, $q\r \mu$ decides
$\phi_\mu^{-1}(i)$ for each $i<n_\mu^p$ (i.e., there are finite sets
$K_{\mu i}^q$ of ordinals, $i<n_\mu^p$, such that
$q\r\mu\forces\phi_\mu^{-1}(i)=K_{\mu i}^q$);
\item"(v)" $q\forces m\in \omega\sm \bigcup_{k<l}
O(T,k)$.
\endroster

Then:  whenever $\rho\in Fn(\o_1\sm M,2)$, 
there is a condition $s\leq 
p,q$ such that $\s^s(\langle\b,m\rangle )=\rho(\b)$ for all $\b\in dom(\rho)$.  
\endproclaim

\demo{Proof} Since $q$ and $p$ are compatible (by hypothesis), we can
choose $r\leq q,p$ such that $\s^{r}$ decides whether or not $n\in A_\b$ 
whenever:
\roster
\item"(1)" $n$ appears in the domain or range of some $\tau_\mu^q$ 
for $\mu\in D^p\cup\{1\}$,
\endroster
and 
\roster
\item"(2)" $\b\in F_1^q\cup(\bigcup\{F_\mu^q\cup K_{\mu i}^q:\mu\in
D^p,i<n_\mu^p\}$.
\endroster

  Let $O=O(\{\tau_\mu^{r}\}_{\mu\in D^p\cup\{1\}},m)$.
 Note that $O$ contains $m$ and is finite (it would 
have to be contained the union of $\{m\}$  and 
the finite set of integers mentioned in the
domain and range of the $\tau_\mu^r$'s). Also $O\cap l=\0$ since 
$r\forces m\not\in  \bigcup_{k<l}O(T,k)$.  

Fix $\rho\in Fn(\o_1\sm M,2)$.  We need to show that there is  
a condition $s\leq 
p,q$ such that $\s^s(\langle\b,m\rangle )=\rho(\b)$ for all $\b\in
dom(\rho)$.  To this end, first 
let  $\s\in Fn(\o_1\times \o,2)$ be such that:
\roster
\item"(3)"$dom(\s)=dom(\s^{r})\cup (dom(\rho)\times O)$;
\item"(4)" $\s(\langle\b,n\rangle )=\rho(\b)$ if $\b\in dom(\rho)$ and $n\in O$; 
\item"(5)" $\s(\langle\b,n\rangle )=\s^{r}(\langle\b,n\rangle )$ otherwise.
\endroster

Now let $$s=\langle\s,\langle\tau_1^q,F_1^q\cup F_1^p\rangle ,\langle\tau_\mu^q,
F_\mu^q\cup F_\mu^p,n_\mu^q\rangle _{\mu\in D^q}\rangle .$$
By (4) above, $s$ (if it is a condition) forces what we want.  
It remains to prove that $s$ is a condition extending both $q$ and
$p$.  

Note that $\s\supset \s^q\cup \s^p$, since:
\roster
\item"(6)"$\s^r\supset \s^q\cup \s^p$;
\item"(7)" $\s=\s^r$ outside of $dom(\rho)\times O$;
\item"(8)"  $dom(\rho)\cap dom(\s^q)=\0$ (recall $q\in M$ and 
$dom(\rho)\subset\o_1\sm M$);
\endroster
and
\roster
\item"(9)" $O\cap l=\0$.
\endroster

So  $s$ can be thought of as $q$ with $\s^q$ extended, the
$F_\mu^q$'s enlarged, but the $\tau_\mu^q$'s left the same.  Also recall
$F_\mu^q\supset F_\mu^p\cap M$  for 
every $\mu\in D^p\cup\{1\}$, so $F_\mu^s\sm F_\mu^q\subset \o_1\sm M$
(this is needed to show $s\r \mu\forces F_\mu^s\subset \o_1\sm
\dot{\delta}_\mu$).  
It easily follows that $s$ is a condition and $s\leq q$.  

It remains to prove that $s\leq p$.  From (i)-(iii) and the definition 
of $s$, we have the necessary containments, and $\tau_\mu^s=\tau_\mu^q$,
so what  we need to show is
that $s\r\mu$ forces:
\roster
\item"(a)" $\forall n\in dom(\tau_\mu^q\sm\tau_\mu^p)
\forall \b\in F_\mu^p[n\in
A_\b\iff \tau_\mu^q(n)\in A_\b];$
\item"(b)" $\forall n\in dom(\tau_\mu^q\sm\tau_\mu^p)
\forall k\leq n_\mu^p[n\in \bigcap_{\phi_\mu (\b)<k}A_\b
\iff \tau^q(n)\in \bigcap_{\phi_\mu (\b)<k}A_\b]$
\endroster
for each $\mu\in D^p$, and 
\roster
\item"(c)"$\forall n\in dom(\tau_1^q\sm\tau_1^p)\forall \b\in F_1^p
[n\in A_\b\iff \tau_1^q(n)\in A_\b\times\o].$
\endroster

 {\it Proof of (a).}  Let $n\in dom(\tau_\mu^q\sm \tau_\mu^p)$, $\b\in
F_\mu^p$, where $\mu\in D^p$.  Since $\tau^q_\mu(n)=\tau_\mu^{r}(n)$,
by the definition of $O$ we have $n\in O\iff \tau_\mu^q(n)\in O$.  
 So if $\langle\b,n\rangle $ is in $ dom(\rho)\times O$, then so is 
$\langle\b,\tau^q_\mu(n)\rangle $ and  by (4) above, $\s (=s\r 1)$ forces $``n\in 
A_\b\iff \rho(\b)=1\iff\tau_\mu^q(n)\in A_\b"$. 
On the other hand, if $\langle\b,n\rangle $ is not in $ dom(\rho)\times O$,
neither is $\langle\b,\tau_\mu^q(n)\rangle $, and so $\s$ and $\s^{r}$ agree on
these two values.  Since $r\leq p$ and $n\in
dom(\tau_\mu^{r}\sm \tau_\mu^p)$, we must have $\s^{r}\forces ``n\in
A_\b\iff \tau_\mu^q(n)\in A_\b"$.  So $\s$, hence $s\r\mu$, forces this too.

  {\it Proof of (c).}  This proof is virtually the same as for (a), 
putting $\mu=1$ and putting $\pi_1(\tau_1^q(n))$ in place of $\tau_\mu^q(n)$.  

 {\it Proof of (b).}  Let $n\in dom(\tau_\mu^q\sm\tau_\mu^p)$ and 
$k\leq n_\mu$.  By (iv), and since 
$r,s\leq q$, $$r\r\mu,s\r\mu\forces``\phi_\mu(\b)<k\iff \b\in 
\bigcup_{i<k}K_{\mu
i}^q".$$  Since $q\in M$ and $dom(\rho)\cap M=\0$, we have 
$dom(\rho)\cap(  \bigcup_{i<k}K_{\mu i}^q)=\0$.  By  the
definition of $r$, $\s^{r}$ decides the values of $\langle\b,n\rangle $  and  
$\langle\b,\tau_\mu^q(n)\rangle $ for any $\b\in\bigcup_{i<k}K_{\mu
i}^q$, and by  definition $\s(=s\r 1)$ decides these values 
the same way.  It follows that $r\r\mu\forces`` n\in 
\bigcap_{\phi_\mu (\b)<k}A_\b" \iff s\r\mu\forces``n\in 
\bigcap_{\phi_\mu (\b)<k}A_\b"$ and 
$r\r\mu\forces`` \tau_\mu^q(n)\in 
\bigcap_{\phi_\mu (\b)<k}A_\b" \iff s\r\mu\forces``\tau_\mu^q(n)\in 
\bigcap_{\phi_\mu (\b)<k}A_\b"$.  Now since $r\r\mu$ forces 
$$``n\in 
\bigcap_{\phi_\mu (\b)<k}A_\b\iff 
 \tau_\mu^q(n)\in\bigcap_{\phi_\mu (\b)<k}A_\b"$$
(since $n\in dom(\tau_\mu^r\sm\tau_\mu^p)$ and $r\leq p$) 
so does $s\r\mu$, and we are done. \qed

\enddemo

The next lemma shows that the hypothesis of Lemma 4.2 is satisfied with
respect to the independent family $\Cal A$ added by the first factor.

\proclaim{Lemma 4.7}  Let $\a<\o_1$.  Suppose $\forces_\a\dot{X}\in \Cal
F^+$.  Then there exists $\g<\o_1$ such that $\forces_\a 
\dot{X}\cap L_\rho\neq \0$ for
every $\rho\in Fn(\o_1\sm\g,2)$ (where $L_\rho=\bigcap_{\b\in dom(\rho)}
A_\b^{\rho(\b)}$).  
\endproclaim
\demo{Proof}  If not, then for each $\g<\o_1$, there is $p_\g\in P_\a$
and $\rho_\g\in Fn(\o_1\sm\g,2)$ such that $p_\g\forces \dot{X}\cap
L_{\rho_\g}=\0$.  Say
$$p_\g=\langle\s^\g,\langle\tau_1^\g,F_1^\g\rangle ,\langle\tau_\mu^\g,F_\mu^\g,n_\mu^\g\rangle _{\mu\in
D^\g}\rangle .$$
We may assume: 
\roster
\item"(i)" The $\s^\g$'s form a delta system with root $\s$;
\item"(ii)" $\exists D\subset \a$ such that $D^\g=D$ for all $\g$;
\item"(iii)" For all $\mu\in D\cup\{1\}$, there is $\tau_\mu$ such that
$\tau_\mu^\g=\tau_\mu$ for all $\g$, and for all $\mu \in D$, there is
$n_\mu$ such that $n_\mu^\g=n_\mu$ for all $\g$;
\item"(iv)" For all $\mu\in D\cup\{1\}$, the $F_\mu^\g$'s form a
$\Delta$-system with root $F_\mu$.
\endroster

Let $p=\langle\s,\langle\tau_1,F_1\rangle ,\langle\tau_\mu,F_\mu,n_\mu\rangle _{\mu\in D}\rangle $.  Let us
write $q\supset p$ if containment holds on each coordinate: i.e.,
$\s^q\supset \s$, $\tau_1^q\supset \tau_1$, etc.  Note that this does
not necessarily imply that $q\leq p$ in the sense of the poset $P_\a$.  

Let $M$ be a countable elementary submodel containing
$p,\dot{X},P_\a$, etc.  Choose $\g\in \o_1\sm M$ such that, for each
$\mu \in D\cup\{1\}$, $$[(F_\mu^\g\sm F_\mu)\cup dom(\s^\g\sm \s)]\cap M=\0.$$
Let $l$ be greater than any integer mentioned by $p_\g$ (i.e., greater
than any $n_\mu$, $\mu\in D\cup\{1\}$, or anything in the domain or range of 
$\s^\g$ or of any $\tau_\mu^\gamma$).  

Let $T=\{t_\mu:\mu\in D\cup\{1\}\}$, where $t_\mu$ is the isomorphism
added by the $\mu^{th}$ coordinate of the forcing.  
Since $\forces_\a\dot{X}\in \Cal F^+$, by Lemma 4.5 there are $r\leq p_\g$  
and $m\in\o$ such that $$r\forces ``m\in \dot{X}\sm
\bigcup_{k<l}O(T,k)".$$  

Consider a maximal antichain of conditions $q\supset p$ such that:
\roster
\item"(a)" $q\forces ``m\in \dot{X}\sm
\bigcup_{k<l}O(T,k)";$
\item"(b)" For each $\mu\in D\cup\{1\}$, $q\r \mu$ decides
$\phi_\mu^{-1}(i)$ for each $i<n_\mu$ (i.e., there are finite sets
$K_{\mu i}^q$ of ordinals, $i<n_\mu$, such that
$q\r\mu\forces"\phi_\mu^{-1}(i)=K_{\mu i}^q$);
\endroster
Any such maximal antichain must contain a condition compatible with $r$.
By $CCC$, and elementarity of $M$, there is such a such a maximal
antichain in $M$, and hence we can choose 
a $q\supset p$ in $M$ satisfying (a) and (b) which is compatible with $r$.  

Note that the conditions of Lemma 4.6 are satisfied for this $q$ with
$p=p_\g$ and $\rho=\rho_\g$.  So there is $s\leq p,q$ with $\s^s(\langle\b,m\rangle )
=\rho_\g(\b)$ for all $\b\in dom(\rho_\g)$.  Then $s\forces m\in
\dot{X}\cap L_{\rho_\g}$.  This contradicts $p_\g\forces \dot{X}\cap
L_{\rho_\g}=\0$, which completes the proof.\qed
\enddemo

At this point, we can see that the iteration can be continued as
claimed, with $\dot{Q_\a}$, $\a> 1$,  forcing $\Cal F$ to be
isomorphic to its restriction to $X_\a$.  Since each $X\in \Cal F^+$ in
the final model is named by some $\dot{X_\a}$, $\a<\o_1$, 
$\Cal F$ is homogeneous in $V^{P_{\o_1}}$.  The forcing $\dot{Q_1}$ 
assures us that $\Cal F$ is isomorphic to $\Cal F\times\o$.

It remains to prove that $\Cal F$ also satisfies condition (iii) in
Theorem 4.1.   Let us note that any subset $I$ of $\o$ which appears
by stage $1<\a<\o_1$ and is in $\Cal F^+$ in $V^{P_\a}$  
remains in $\Cal F^+$ in $V^{P_{\o_1}}$ (since for each finite $H\subset 
\o_1$, the statement ``$I\cap \bigcap_{\b\in H}A_\b\neq \0$" is absolute
for transitive models of set theory).

 Therefore, since any $f:\o\to \o$ in the final model appears in some 
initial stage,  it suffices to prove that condition (iii) holds
at every stage $1<\a<\o_1$.  So suppose $$\forces_\a ``\dot{f}:\o\to \o
\text{ and } \forall F\in \Cal F(\dot{f}(F)\text{ is unbounded })".$$ 
Let $M$ be a countable elementary submodel containing $P_\a,
\dot{f},...$.  Let $M\cap \o_1=\{\a_0,\a_1,...\}$.  

Let $T_k=\{t_{\a_l}:l<k, 1\leq \a_l<\a\}$.  Let $a_\0=\0$.  
By induction on $k\geq 1$, we can define, for each $\s\in \o^k$, conditions 
$a_\s\in M\cap P_\a$, and $i_\s,j_\s\in \o$ such that:
\roster 
\item"(a)" $a_\s\forces ``\dot{f}(i_\s)=j_\s\text{ and }i_\s\in
\bigcap_{l<|\s|}A_{\a_l}\sm \bigcup_{m<|\s|}O(T_{|\s|},m)"$;
\item"(b)" $i_{\s^\frown  \langle n\rangle }>i_\s$ and $j_{\s^\frown  \langle n\rangle }>j_\s$;
\item"(c)" $D^{a_\s}\supset \{\a_l:l<|\s|,1<\a_l<\a\}$; 
\item"(d)" For each $\mu\in D^{a_\s}$, $n_\mu^{a_\s}\geq |\s|$, and there is
$\d_\mu^{a_\s}<\o_1$ such that $a_\s\r\mu\forces
\dot{\d_\mu}=\d_\mu^{\a_\s}$;
\item"(e)" For each $\mu\in D^{a_\s}$, for each $i<|\s|$, 
$a_\s\r \mu$ decides the value of $\phi_\mu^{-1}(i)$;
\item"(f)" For each $\mu\in D^{a_\s}\cup \{1\}$,
$dom(\tau_\mu^{a_\s}\supset |\s|$;
\item"(g)" $F_1^{a_\s}\supset \{\a_l:l<|\s|\}$ and for each $\mu\in 
D^{\a_\s}$, $F_\mu^{a_\s}\supset \{\a_l:l<|\s|, \d_\mu^{a_\s}\leq
\a_l<M\cap \o_1\}$; 
\item"(h)" $a_\s\leq a_{\s\r l}$ for each $ l\leq |\s|$;
\item"(j)" $\{a_{\s^\frown  \langle n\rangle }:n<\o\}$ is a maximal antichain below
$a_\s$.
\endroster

Note that (a) and (b) can be obtained  because $\dot{f}$ is forced to be 
unbounded on every $F\in \Cal F$.  Then (c)-(g) can be obtained simply 
by extending as necessary.   So given that $a_\s$ has been defined, there
is a maximal antichain of conditions below $a_\s$ satisfying (a)-(g), 
and by $CCC$, there is one in $M$; hence (h) and (j) can be obtained as
well.  

Note that for any $P_\a$-generic $G$, there is 
a unique $h\in\o^\o$ with $a_{h\r n}\in G$ for each $n\in \o$.  So there  
a corresponding set $I=\{i_{h\r n}:n<\o\}$ in $V[G]$.  By (a) and (b), 
$\dot{f}_G$ is strictly increasing on $I$.  

Let $\dot{I}$ be
a $P_\a$-name for $I$.  It remains to prove that $\forces_\a\dot{I}\in
\Cal F^+$.  Suppose not.  Then there is some $p\in P_\a$ and finite
subset $H$ of $\o_1$ such that $p\forces ``\dot{I}\cap\bigcap_{\b\in
H}A_\b=\0"$.

Note that there is some $a_\s$ compatible with $p$ for $\s$ of any
prescribed length.  So we can choose $a_\s$ compatible with $p$ so that:
\roster
  \item"(k)" $|\s|$ is greater than any integer mentioned in $p$;
\item"(l)" $\{\a_l:l<|\s|\}\supset D^p\cup (H\cap M)\cup (\cup\{F_\mu^p\cap M:
\mu\in D^p\})$;
\endroster

Then from (l) and (a) we have:
\roster
\item"(m)" $a_\s\forces i_\s\in \bigcap_{\g\in H\cap M}A_\g$.
\endroster
By (c) and (l) we have:
\roster
\item"(n)" $D^{a_\s}\supset D^p$.
\endroster
By (k), (n), (f), and compatibility we have:
\roster
\item"(o)" For each $\mu\in D^p\cup\{1\}$,  $\tau_\mu^{a_\s}\supset
\tau_\mu^p$.
\endroster
By (d), (g), (l), and compatibility we have:
\roster
\item"(p)" For each $\mu\in D^p\cup\{1\}$, $F_\mu^{a_\s}\supset
F_\mu^p\cap M$, and $n_\mu^{a_\s}\geq n_\mu^p$ for each $\mu\in D^p$.
\endroster

It follows from (a),(e),(k),(n),(o), and (p) that the conditions of Lemma 4.6
are satisfied for this $p$ with $q=a^\s$, $m=i_\s$, and $l=|\s|$.  So,
there exists $s\leq a_\s,p$ with $\s^s(\langle \b,i_\s\rangle )=1$ for each $\b\in
H\sm M$.  Then by (m) and $s\leq a_\s$, we have $$s\forces 
i_\s\in \dot{I}\cap \bigcap_{\b\in H}A_\b.$$  This contradicts $s\leq
p$, and completes the proof of Theorem 4.1.  \qed

We wish to remark here that there is a theorem of W. H. Woodin 
\cite{W} which suggests
that the Continuum Hypothesis implies that a filter with the properties
given in Theorem 4.1 exists.
In particular he has shown that if $\phi$ is a $\Sigma^2_1$ sentence
(i.e. a sentence of the form
$\exists \Cal F \subseteq P(\omega) Q_1 x_1 \ldots Q_k x_k 
\psi(\Cal F,x_1,\ldots,x_k)$
where $Q_i$ is a quantifier over the reals and $\psi$ is a
quantifier free formula)
then, provided that certain large cardinals exist
(a supercompact is sufficient),  if $\phi$ can be forced to be true,
then it holds in any forcing extension satisfying the Continuum Hypothesis.
It is not difficult to show that the existence of a filter
with the properties of Theorem 4.1 is a statement
of this form.
If the existence of such a filter were not a consequence
of the Continuum Hypothesis then this would be of independent interest.

\head 5. Concluding remarks. \endhead

  The topology of our spaces $T_\a$ are determined by the filter $\Cal F$ 
on $\o$,  and by the ``$\cup_{\Cal F}$ construction" for defining the neighborhoods
of the point $\infty$ which we added to the countably many 
copies $\{n\}\times T_\a$, $n<\o$, of $T_\a$, in order to construct $T_{\a+1}$ 
from $T_\a$.   A different natural way of defining the topology 
of $T_{\a+1}$ is to declare $N$
to be a neighborhood of $\infty$ iff $\infty\in N$ and $N$ contains a
neighborhood of the point $(n,\infty)$ of maximal rank in $\{n\}\times
T_\a$, for $\Cal F$-many $n\in \o$.  Let $T_\a'$ be the spaces obtained
by defining them inductively this way, i.e., just like the $T_\a$'s but
using this finer neighborhood base for $\infty$.  Then $T_2'$ is the
same as $T_2$, i.e., it's homeomorphic to the space $\o\cup \{\Cal F\}$ with $\o$ the
set of isolated points and a neighborhood of $\Cal F$ having the form
$F\cup\{\Cal F\}$, $F\in \Cal F$.  But $T_3'$ is different; in
particular, the restriction of the neighborhood filter of $\infty$ to the 
isolated points is isomorphic to the filter $\Cal F^2$ on $\o^2$ defined by 
$$\Cal F^2=\{A\subset \o^2:\{n:\{m:(n,m)\in A\}\in\Cal F\}\in \Cal
F\}.$$  A filter $\Cal F$ is said to be {\it idempotent} if it is isomorphic to 
$\Cal F^2$.  

The following result shows that, while the $T_\a'$ construction does yield a
3-Toronto space in $ZFC$ (iff $\Cal F$ is an ultrafilter), $T_\a'$ is
{\it never} an $\a$-Toronto space if $\a\geq 4$.  

\proclaim{Theorem 5.1} Let $\Cal F$ be a filter on $\o$, and let 
$T_\a'$, $\a$ an ordinal, be the spaces defined as above.  Then:
\roster
\item"(a)" $T_3'$ is a 3-Toronto space iff $\Cal F$ is an ultrafilter;
\item"(b)" $T_\a'$ is not an $\a$-Toronto space if $\a\geq 4$.
\endroster
\endproclaim
\demo{Proof} (a) We can consider $T_3'$ as the set $(\o\times (\o+1))\cup
\{\infty\}$, where $\o\times \o$ is the set of isolated points, each
$\<n,\o\>$ has a neighborhood base of sets of the form $\{n\}\times
(F\cup\{\o\})$,
$F\in \Cal F$, and a neighborhood of $\infty$
consists of $\infty$ together with neighborhoods of $\<n,\o\>$ for $\Cal
F$-many $n$.  

Suppose $\Cal F$ is an ultrafilter and $X\subset T_3'$ has rank 3.  Then
for $\Cal F$-many $n$, it must be the case that there are
$F_n\in\Cal F$ such that $\{n\}\times (F_n\cup\{\o\})\subset X$.  Now,
using the easily verified fact that any ultrafilter is homogeneous, it
is not difficult to prove that $X\cong T_3'$.  

Suppose on the other hand that $\Cal F$ is not an ultrafilter. 
Let $A\subset \o$ be such that
$A,\o\sm A\in \Cal F^+$.  Let $X=T_3'\sm (A\times \{\o\})$.  Then $X$ has rank
3.  But the subset $A\times \o$ of $X$ is a set of isolated points whose
 only limit point in $X$ is $\infty$, while in $T_3'$ any set of isolated
points clustering at $\infty$ must also have cluster points of the form
$ \langle n,\o\rangle $.  So $X$ is not homeomorphic to $T_3'$.  

(b) Let $\a\geq 4$, and consider removing the level one points from a 
neighborhood of some level 2 point of $T_\a'$.  By the discussion preceding the 
theorem, in the resulting subspace this level 2 point of $T_\a'$ becomes a level 
one point with a neighborhood filter isomorphic to $\Cal F^2$.  
Thus it must be the case that $\Cal F$ is an idempotent filter.  
Frol\'ik showed in \cite{F} that no ultrafilter is idempotent, 
so $\Cal F$ is not an ultrafilter.
Note that any set of isolated points of $T_\a'$ which has a limit point
at level 2 must also have a limit point at level 1.  But, using an idea
from the proof of part (a), we see that there is a subspace of $T_\a'$
of rank $\a$ which fails to have this property. Hence $T_\a'$ is not 
$\a$-Toronto.  \qed
\enddemo

{\bf Remark.} An idempotent filter on $\o$ has been constructed in $ZFC$ 
by M. Katetov
\cite{Ka}, and a homogeneous one (also in $ZFC$) by J. Steprans (unpublished
note).  The latter result answers a question of Steprans in \cite{S}.

\enddocument